\documentclass[english]{article}
\usepackage[T1]{fontenc}
\usepackage[latin9]{inputenc}
\usepackage{geometry}
\usepackage{amsmath}
\usepackage{amssymb}
\usepackage{esint}

\makeatletter
\@ifundefined{definecolor}
 {\usepackage{color}}{}
\@ifundefined{definecolor}{\usepackage{color}}{}

\usepackage{babel}
\usepackage{nicefrac}
\usepackage{sagetex}

\makeatother

\usepackage{babel}
\begin{document}

\title{Computing Spectral Elimination Ideals}

\author{Edinah K. Gnang%
\thanks{School of Mathematics, Institute for Advanced Study. Email: gnang@ias.edu.%
}}
\maketitle
\begin{abstract}
We present here an overview of the hypermatrix spectral decomposition
deduced from the Bhattacharya-Mesner hypermatrix algebra \cite{BM1,BM2}.
We describe necessary and sufficient conditions for the existence
of a spectral decomposition. We further extend to hypermatrices the
notion of resolution of identity and use them to derive hypermatrix
analog of matrix spectral bounds. Finally we describe an algorithm
for computing generators of the spectral elimination ideals which
considerably improves on Groebner basis computation suggested in \cite{GER}.
\end{abstract}

\section{Introduction}

This brief note discusses \emph{hypermatrices}, a generalization of
matrices which corresponds to a finite set of numbers each of which
is indexed by a member of an integer cartesian product set of the
form $\left\{ 0,\cdots,\left(n_{0}-1\right)\right\} \times\cdots\times\left\{ 0,\cdots,\left(n_{l-1}-1\right)\right\} $.
Such a hypermatrix is said to be of order $l$ and more conveniently
called an $l$-hypermatrix. The algebra and the spectral analysis
of hypermatrices arise from generalizations of familiar concepts of
linear algebra. The reader is referred to \cite{L} for a survey of
important hypermatrix results. The reader is also refereed to \cite{LQ}
for a detailed survey of the various approaches to the spectral analysis
of hypermatrices. The algebra discussed here differs considerably
from the hypermatrix algebras surveyed in \cite{L}. The hypermatrix
algebra discussed here centers around the Bhattacharya-Mesner hypermatrix
product operation motivated by generalizations of association schemes
introduced in \cite{BM1,BM2,B} and followed up in \cite{GER}. Although
the scope of the Bhattacharya-Mesner algebra extends to hypermatrices
of all integral orders, the present discussion will be restricted
for notational convenience to $3$-hypermatrices since all the result
presented here generalize straight-forwardly to greater order hypermatrices.\\
Our main result provides necessary and sufficient conditions for the
existence of a spectral decomposition for a given hypermatrix which
is symmetric under cyclic permutation of it's indices. We also describe
how to extend to hypermatrices the notion of resolution of identity
deduced from orthogonal hypermatrices introduced in \cite{GER}. We
also extend to hypermatrices the symmetrization approach to Singular
Value Decomposition. We further derive hypermatrix analog of matrix
spectral bounds. Finally, we describe an algorithm for computing generators
of the elimination ideals which considerably improves on Groebner
basis computation approach suggested in \cite{GER} for computing
elimination ideals.

\section{Hypermatrix orthogonality}

We describe here a spectral decomposition for 3-hypermatrices, deduced
from the Bhattacharya-Mesner algebra. The proposed spectral decomposition
builds on the notion of hypermatrix orthogonality introduced in \cite{GER},
defined for arbitrary order hypermatrices. In particular, we establish
the existence of arbitrary order orthogonal hypermatrices by describing
an explicit parametrization orthogonal hypermatrices resulting from
direct sums of hyermatrices of size $2\times2\times\cdots\times2\times2$
( the direct sum refers here to hypermatrix diagonal block construction
quite analogous to the matrix counterpart). We recall that hypermatrix
orthogonality for an $m$-hypermatrix is determined by the constraints
\begin{equation}
\boldsymbol{\Delta}=\bigcirc_{0\le t<m}\left(\mathbf{Q}^{T^{(m-t)}}\right)
\end{equation}
which is more explicitly expressed as 
\begin{equation}
\delta_{i_{0}\cdots i_{m-1}}=\sum_{0\le k<n}q_{i_{0}ki_{2}\cdots i_{m-2}i_{m-1}}\cdot q_{i_{1}ki_{3}\cdots i_{m-1}i_{0}}\cdot q_{i_{2}ki_{4}\cdots i_{0}i_{1}}\cdot\:\cdots\:\cdot q_{i_{m-2}ki_{0}\cdots i_{m-4}i_{m-3}}\cdot q_{i_{m-1}ki_{1}\cdots i_{m-3}i_{m-2}},
\end{equation}
(where $\delta_{i_{0}\cdots i_{m-1}}$ denotes the entries of the
Kronecker delta). It follows that the parametrization of the sought
after family of direct sums of orthogonal hypermatrices is completely
determined by the parametrization of hypermatrices of dimensions $\underbrace{2\times2\times\cdots\times2\times2}_{m\mbox{ operands}}$
which we determine by solving the linear constraints of the form
\[
\ln q_{i_{0}0i_{2}\cdots i_{m-2}i_{m-1}}+\ln q_{i_{1}0i_{3}\cdots i_{m-1}i_{0}}+\ln q_{i_{2}0i_{4}\cdots i_{0}i_{1}}+\:\cdots\:+\ln q_{i_{m-2}0i_{0}\cdots i_{m-4}i_{m-3}}+\ln q_{i_{m-1}0i_{1}\cdots i_{m-3}i_{m-2}}=
\]
\begin{equation}
i\pi+\ln q_{i_{0}1i_{2}\cdots i_{m-2}i_{m-1}}\ln q_{i_{1}1i_{3}\cdots i_{m-1}i_{0}}+\ln q_{i_{2}1i_{4}\cdots i_{0}i_{1}}+\:\cdots\:+\ln q_{i_{m-2}1i_{0}\cdots i_{m-4}i_{m-3}}+\ln q_{i_{m-1}1i_{1}\cdots i_{m-3}i_{m-2}}.
\end{equation}
Note that there will be one constraint for every orbit of the action
of the cyclic group on $m$-tuples. For instance the constraints above
yield the following parametrization for orthogonal hypermatrices of
size $2\times2\times2$ expressed by 
\begin{equation}
q_{000}=\frac{e^{r_{3}}}{\left(e^{\left(3\, r_{3}\right)}+e^{\left(3\, r_{6}\right)}\right)^{\frac{1}{3}}},\, q_{001}=e^{r_{4}},\, q_{010}=\frac{e^{r_{6}}}{\left(e^{\left(3\, r_{3}\right)}+e^{\left(3\, r_{6}\right)}\right)^{\frac{1}{3}}},\, q_{011}=e^{r_{2}}
\end{equation}
\begin{equation}
q_{100}=-e^{\left(r_{2}-r_{3}-r_{4}+r_{5}+r_{6}\right)},\, q_{101}=\frac{e^{\left(r_{1}+r_{3}-r_{6}\right)}}{\left(e^{\left(3\, r_{1}\right)}+e^{\left(3\, r_{1}+3\, r_{3}-3\, r_{6}\right)}\right)^{\frac{1}{3}}},\, q_{110}=e^{r_{5}},\, q_{111}=\frac{e^{r_{1}}}{\left(e^{\left(3\, r_{1}\right)}+e^{\left(3\, r_{1}+3\, r_{3}-3\, r_{6}\right)}\right)^{\frac{1}{3}}}
\end{equation}
for arbitrary choice of values of parameters $\left\{ r_{k}\right\} _{0<k<7}$.

\section{Matrix spectral elimination ideals}

We recall that the spectral constraint for symmetric real matrices
are expressed as follows 
\begin{equation}
\begin{cases}
\begin{array}{ccc}
\mathbf{A} & = & \left(\mathbf{Q}\cdot\mathbf{D}\right)\cdot\left(\mathbf{Q}\cdot\mathbf{D}\right)^{T}\\
\left[\mathbf{Q}\cdot\mathbf{Q}^{T}\right]_{i,j} & = & \begin{cases}
\begin{array}{cc}
1 & \mbox{if }i=j\\
0 & \mbox{otherwise }
\end{array} & \forall\:0\le i,j<n\end{cases}\\
\mathbf{D}^{\star^{2}} & = & \mathbf{D}^{T}\cdot\mathbf{D}
\end{array}\end{cases},
\end{equation}
and the corresponding invariance formulation is given by 
\begin{equation}
\mathbf{A}\cdot\left[\left(\mathbf{Q}\cdot\mathbf{D}\right)^{T}\right]^{-1}=\mathbf{Q}\cdot\mathbf{D},
\end{equation}
provided of course that the matrix $\left(\mathbf{Q}\cdot\mathbf{D}\right)^{T}$
is invertible. We describe here the determination of generators for
the elimination ideals
\begin{equation}
I_{\mathbf{D}}=\left\rangle \mathbf{A}-\left(\mathbf{Q}\cdot\mathbf{D}\right)\cdot\left(\mathbf{Q}\cdot\mathbf{D}\right)^{T},\,\mathbf{Q}\cdot\mathbf{Q}^{T}-\boldsymbol{\Delta},\,\mathbf{D}^{\star^{2}}-\mathbf{D}^{T}\cdot\mathbf{D}\right\langle \cap\mathbb{C}\left[\mathbf{D}\right]
\end{equation}
 and 
\begin{equation}
I_{\mathbf{Q}}=\left\rangle \mathbf{A}-\left(\mathbf{Q}\cdot\mathbf{D}\right)\cdot\left(\mathbf{Q}\cdot\mathbf{D}\right)^{T},\,\mathbf{Q}\cdot\mathbf{Q}^{T}-\boldsymbol{\Delta},\,\mathbf{D}^{\star^{2}}-\mathbf{D}^{T}\cdot\mathbf{D}\right\langle \cap\mathbb{C}\left[\mathbf{Q}\right]
\end{equation}
where $\left\rangle \mathbf{A}-\left(\mathbf{Q}\cdot\mathbf{D}\right)\cdot\left(\mathbf{Q}\cdot\mathbf{D}\right)^{T},\mathbf{Q}\cdot\mathbf{Q}^{T}-\boldsymbol{\Delta},\mathbf{D}^{\star^{2}}-\mathbf{D}^{T}\cdot\mathbf{D}\right\langle $
denotes the ideal generated by the matrix spectral constraints. While
the characterization of both these elimination ideals is well known
in the case of matrices, our aim is to present the derivation of these
elimination ideal so to suggest a natural generalization of the derivation
to hypermatrices of all integral orders. We emphasize that the proposed
derivation avoids the Groebner basis computations suggested in \cite{GER}.

As a starting point for the derivation we consider the following equivalent
formulation of the matrix spectral constraints 
\begin{equation}
\begin{cases}
\begin{array}{ccc}
a_{ij} & = & \left\langle \left(\mathbf{q}_{i}\star\boldsymbol{\lambda}\right),\left(\boldsymbol{\lambda}\star\mathbf{q}_{j}\right)\right\rangle \\
\delta_{ij} & = & \left\langle \mathbf{q}_{i},\mathbf{q}_{j}\right\rangle 
\end{array} & \forall\:0\le i\le j<n\end{cases}
\end{equation}
where $\boldsymbol{\lambda}$ denotes the vector whose entries are
the principal square roots of the eigenvalues of $\mathbf{A}$. We
therefore deduce the following expressions for the columns of $\mathbf{A}$
\begin{equation}
\mathbf{a}_{j}=\left(a_{ij}=\left\langle \left(\mathbf{q}_{i}\star\boldsymbol{\lambda}\right),\left(\boldsymbol{\lambda}\star\mathbf{q}_{j}\right)\right\rangle \right)_{0\le i<n}\Rightarrow\mathbf{a}_{j}=\mathbf{Q}^{T}\cdot\left(\boldsymbol{\lambda}^{\star^{2}}\star\mathbf{q}_{j}\right)
\end{equation}
and in particular 
\begin{equation}
\left\langle \mathbf{a}_{i},\mathbf{a}_{j}\right\rangle =\left[\mathbf{Q}^{T}\cdot\left(\boldsymbol{\lambda}^{\star^{2}}\star\mathbf{q}_{j}\right)\right]^{T}\cdot\left[\mathbf{Q}^{T}\cdot\left(\boldsymbol{\lambda}^{\star^{2}}\star\mathbf{q}_{i}\right)\right]
\end{equation}
\begin{equation}
\Rightarrow\left\langle \mathbf{a}_{i},\mathbf{a}_{j}\right\rangle =\left\langle \left(\boldsymbol{\lambda}^{\star^{2}}\star\mathbf{q}_{j}\right),\left(\boldsymbol{\lambda}^{\star^{2}}\star\mathbf{q}_{i}\right)\right\rangle 
\end{equation}
quite similarly let $\mathbf{a}_{j}^{[k-1]}$ denote the $j$-th column
of the matrix power $\mathbf{A}^{k-1}$, we have 
\begin{equation}
\left\langle \mathbf{a}_{i},\mathbf{a}_{j}^{[k-1]}\right\rangle =\left\langle \left(\boldsymbol{\lambda}^{\star^{k}}\star\mathbf{q}_{j}\right),\left(\boldsymbol{\lambda}^{\star^{k}}\star\mathbf{q}_{i}\right)\right\rangle =\left\langle \boldsymbol{\lambda}^{\star^{2k}},\left(\mathbf{q}_{j}\star\mathbf{q}_{i}\right)\right\rangle .
\end{equation}
We therefore deduce Vandermonde block linear constraints which we
conveniently express in matrix form as follows 
\begin{equation}
\left(\begin{array}{c}
\left[\mathbf{q}_{i}\star\mathbf{q}_{j}\right]_{0}\\
\vdots\\
\left[\mathbf{q}_{i}\star\mathbf{q}_{j}\right]_{n-1}
\end{array}\right)_{0\le i\le j<n}=\left(\mathbf{I}_{{n+1 \choose 2}}\otimes\mathbf{V}\left(\boldsymbol{\lambda}^{\star^{2}}\right)\right)^{-1}\cdot\left(\begin{array}{c}
\left[\mathbf{A}^{0}\right]_{i,j}\\
\vdots\\
\left[\mathbf{A}^{n-1}\right]_{i,j}
\end{array}\right)_{0\le i\le j<n}
\end{equation}
where $\mathbf{V}\left(\mathbf{x}\right)$ denotes the $n\times n$
Vandermonde matrix expressed by 
\begin{equation}
\mathbf{V}\left(\mathbf{x}\right):=\left(v_{ij}\left(\mathbf{x}\right)=\left(x_{j}\right)^{i}\right)_{0\le i,j<n}.
\end{equation}
Finally, the elimination ideal $I_{\mathbf{D}}$ is determined by
equating corresponding expressions to obtain ${n \choose 2}$ vector
constraints determined by the equalities 
\begin{equation}
\left\{ \left(\mathbf{q}_{i}\star\mathbf{q}_{j}\right)^{\star^{2}}=\mathbf{q}_{i}^{\star^{2}}\star\mathbf{q}_{j}^{\star^{2}}\right\} _{0\le i<j<n}
\end{equation}
We have therefore derived generators for the ideal of elementary symmetric
polynomials in the square roots of the eigenvalues of $\mathbf{A}$. 

Although the derivation steps described in the previous paragraph
are insightful for the matrix case, unfortunately, these derivation
steps are of limited interest for higher order hypermatrices. In fact
the derivation steps described here only extend to hypermatrices which
correspond to direct sums of hypermatrices whose size is of the form
$2\times2\times2\times\cdots\times2$. Fortunately, however, the derivation
of the elimination ideal $I_{\mathbf{Q}}$ naturally extend to general
hypermatrices. As a result, we advocate the use of the elimination
ideal $I_{\mathbf{Q}}$, as a basis for iterative procedures for approximating
the spectral decomposition of hypermatrices. Our starting point for
the matrix case will be $2{n+1 \choose 2}$ quadratic constraints
\begin{equation}
\begin{cases}
\begin{array}{ccc}
a_{ij} & = & \left\langle \left(\mathbf{q}_{i}\star\boldsymbol{\lambda}\right),\left(\boldsymbol{\lambda}\star\mathbf{q}_{j}\right)\right\rangle \\
\delta_{ij} & = & \left\langle \mathbf{q}_{i},\mathbf{q}_{j}\right\rangle 
\end{array} & \forall\:0\le i\le j<n\end{cases}.
\end{equation}
The main step of the derivation consists in combining the decomposition
constraints with the orthogonality constraints via the use of induced
resolutions of identity. We recall for the convenience of the reader
that the resolution of identity induced by the orthogonal matrix $\mathbf{Q}$
can be expressed as follows 
\begin{equation}
\forall\:\mathbf{u},\mathbf{v}\in\mathbb{C}^{n},\qquad\left\langle \mathbf{u},\mathbf{v}\right\rangle =\left\langle \mathbf{u},\mathbf{v}\right\rangle _{\left(\sum_{0\le t<n}\mathbf{q}_{t}\cdot\mathbf{q}_{t}^{T}\right)}=\sum_{0\le t<n}\left\langle \mathbf{u},\mathbf{v}\right\rangle _{\left(\mathbf{q}_{t}\cdot\mathbf{q}_{t}^{T}\right)}.
\end{equation}
The resolution of identity property is precisely the reason why orthogonality
plays such a crucial role in the formulation of the spectral constraints.
Using the resolution of identity, we reduce the spectral constraints
to the ${n+1 \choose 2}$ constraints 
\begin{equation}
\left\{ a_{ij}=\sum_{0\le k<n}\left\langle \left(\mathbf{q}_{i}\star\boldsymbol{\lambda}\right),\left(\boldsymbol{\lambda}\star\mathbf{q}_{j}\right)\right\rangle _{\mathbf{q}_{k}\cdot\mathbf{q}_{k}^{T}}\right\} _{0\le i\le j<n}
\end{equation}
which may more conveniently be rewritten as 
\begin{equation}
\left\{ a_{ij}=\sum_{0\le k<n}\left\langle \left(\boldsymbol{\lambda}\cdot\boldsymbol{\lambda}^{T}\right),\,\left(\mathbf{q}_{k}\cdot\mathbf{q}_{k}^{T}\right)\star\left(\mathbf{q}_{i}\cdot\mathbf{q}_{j}^{T}\right)\right\rangle \right\} _{0\le i\le j<n}.
\end{equation}
We may think off the set of constraints above as ${n+1 \choose 2}$
linear constraints in the ${n+1 \choose 2}$ variables $\left\{ \lambda_{i}\lambda_{j}\right\} _{0\le i\le j<n}$.
The system is thus solved via Cramer's rule and thus the elimination
ideal $I_{\mathbf{Q}}$ is obtain from equating the appropriate constraints,
suggested by the equality
\begin{equation}
\left\{ \left(\lambda_{i}\,\lambda_{j}\right)^{2}=\lambda_{i}^{2}\,\lambda_{j}^{2}\right\} _{0\le i<j<n}.
\end{equation}
Clearly the ideal $I_{\mathbf{Q}}$ is equivalently characterized
by 
\begin{equation}
\forall\;0\le i<j<n,\quad\left[\mathbf{Q}\cdot\mathbf{A}\cdot\mathbf{Q}^{T}\right]_{i,j}=0
\end{equation}
subject to the constraints 
\begin{equation}
\mathbf{Q}^{T}\cdot\mathbf{Q}=\mathbf{I}
\end{equation}
However the advantage of the proposed derivation of the elimination
ideal $I_{\mathbf{Q}}$ is that it completely determines the expression
of the eigenvalues in terms of the eigenvectors and the generators
for the constraint are half the size of the original constraints.
We may further remark that the derivation suggest a natural mapping
between the orthogonality constraints 
\[
\left\{ \left\langle \mathbf{q}_{i},\mathbf{q}_{j}\right\rangle =0\right\} _{0\le i<j<n}
\]
and the constraints which determines $I_{\mathbf{Q}}$, namely 
\begin{equation}
\left\{ \left(\lambda_{i}\,\lambda_{j}\right)^{2}-\lambda_{i}^{2}\,\lambda_{j}^{2}=0\right\} _{0\le i<j<n}
\end{equation}
It is not unlikely that such a mapping may in off itself suggest alternative
proof of existence and unicity of the spectral decomposition for symmetric
matrices.

\section{3-Hypermatrix elimination ideals.}

We describe here how we extend to hypermatrices the elimination ideal
computations describe in the previous section. We recall here that
third order hypermatrix spectral constraints introduced in \cite{GER}
are expressed as 
\begin{equation}
\begin{cases}
\begin{array}{ccc}
\mathbf{A} & = & \circ\left(\circ\left(\mathbf{Q},\mathbf{D},\mathbf{D}^{T}\right),\,\circ\left(\mathbf{Q},\mathbf{D},\mathbf{D}^{T}\right)^{T^{2}},\,\circ\left(\mathbf{Q},\mathbf{D},\mathbf{D}^{T}\right)^{T}\right)\\
\left[\circ\left(\mathbf{Q},\,\mathbf{Q}^{T^{2}},\,\mathbf{Q}^{T}\right)\right]_{i,j,k} & = & \begin{cases}
\begin{array}{cc}
1 & \mbox{if }i=j=k\\
0 & \mbox{otherwise }
\end{array} & \forall\:0\le i,j,k<n\end{cases}\\
\mathbf{D}^{\star^{3}} & = & \circ\left(\mathbf{D}^{T},\,\mathbf{D}^{T^{2}},\,\mathbf{D}\right)
\end{array}\end{cases}
\end{equation}
Just as we did for matrices, we may express for hypermartrices the
corresponding invariance equality expressed by 
\begin{equation}
\circ\left(\mathbf{A},\,\left[\circ\left(\mathbf{Q},\mathbf{D},\mathbf{D}^{T}\right)^{T^{2}}\right]^{\left(-1\right)_{1}},\,\left[\circ\left(\mathbf{Q},\mathbf{D},\mathbf{D}^{T}\right)^{T}\right]^{\left(-1\right)_{2}}\right)=\circ\left(\mathbf{Q},\mathbf{D},\mathbf{D}^{T}\right)
\end{equation}
provided of course that the pair of hypermatrices $\left(\:\circ\left(\mathbf{Q},\mathbf{D},\mathbf{D}^{T}\right)^{T^{2}},\;\circ\left(\mathbf{Q},\mathbf{D},\mathbf{D}^{T}\right)^{T}\right)$
forms an invertible pair in the sense defined in \cite{BM2}. As pointed
out in the previous section, the approach for deriving the elimination
ideal 
\begin{equation}
I_{\mathbf{D}}=\left\rangle \mathbf{A}-\circ\left(\circ\left(\mathbf{Q},\mathbf{D},\mathbf{D}^{T}\right),\circ\left(\mathbf{Q},\mathbf{D},\mathbf{D}^{T}\right)^{T^{2}},\circ\left(\mathbf{Q},\mathbf{D},\mathbf{D}^{T}\right)^{T}\right),\boldsymbol{\Delta}-\circ\left(\mathbf{Q},\mathbf{Q}^{T^{2}},\mathbf{Q}^{T}\right),\mathbf{D}^{\star^{3}}-\circ\left(\mathbf{D}^{T},\mathbf{D}^{T^{2}},\mathbf{D}\right)\right\langle \cap\mathbb{C}\left[\mathbf{D}\right]
\end{equation}
only extends to hypermatrices which are direct sums of hypermatrices
whos size is of the form $2\times2\times2\times\cdots\times2$. In
the particular case of 3-hypermatrices the problem completely reduces
to the spectral decomposition of $2\times2\times2$ hypermatrices
determined by the constraints 
\begin{equation}
\begin{cases}
\begin{array}{ccc}
\left\langle \left(\mathbf{w}_{0}^{\star^{6}}\right),\,\left(\mathbf{q}_{00}\star\mathbf{q}_{00}\star\mathbf{q}_{00}\right)\right\rangle  & = & a_{000}\\
\left\langle \left(\mathbf{w}_{1}^{\star^{6}}\right),\,\left(\mathbf{q}_{11}\star\mathbf{q}_{11}\star\mathbf{q}_{11}\right)\right\rangle  & = & a_{111}\\
\left\langle \left(\mathbf{w}_{0}^{\star^{2}}\star\mathbf{w}_{1}^{\star^{4}}\right),\,\left(\mathbf{q}_{01}\star\mathbf{q}_{10}\star\mathbf{q}_{11}\right)\right\rangle  & = & a_{011}\\
\left\langle \left(\mathbf{w}_{0}^{\star^{4}}\star\mathbf{w}_{1}^{\star^{2}}\right),\,\left(\mathbf{q}_{10}\star\mathbf{q}_{01}\star\mathbf{q}_{00}\right)\right\rangle  & = & a_{100}\\
\left\langle \left(\mathbf{w}_{0}^{\star^{6}}\right)^{\star^{0}},\:\left(\mathbf{q}_{00}\star\mathbf{q}_{00}\star\mathbf{q}_{00}\right)\right\rangle  & = & \delta_{000}\\
\left\langle \left(\mathbf{w}_{1}^{\star^{6}}\right)^{\star^{0}},\:\left(\mathbf{q}_{11}\star\mathbf{q}_{11}\star\mathbf{q}_{11}\right)\right\rangle  & = & \delta_{111}\\
\left\langle \left(\mathbf{w}_{0}^{\star^{2}}\star\mathbf{w}_{1}^{\star^{4}}\right)^{\star^{0}},\:\left(\mathbf{q}_{01}\star\mathbf{q}_{10}\star\mathbf{q}_{11}\right)\right\rangle  & = & \delta_{011}\\
\left\langle \left(\mathbf{w}_{0}^{\star^{4}}\star\mathbf{w}_{1}^{\star^{2}}\right)^{\star^{0}},\:\left(\mathbf{q}_{10}\star\mathbf{q}_{01}\star\mathbf{q}_{00}\right)\right\rangle  & = & \delta_{100}
\end{array}\end{cases}.
\end{equation}
The system of equations above corresponds to a block Vandermonde set
of linear constraints which yields the following constraints 
\begin{equation}
\begin{cases}
\begin{array}{ccc}
\frac{\left(w_{00}^{4}w_{01}^{2}-w_{01}^{4}w_{11}^{2}\right)^{3}\left(a_{000}-w_{01}^{6}\right)}{a_{001}^{3}\left(w_{00}^{6}-w_{01}^{6}\right)} & = & \frac{\left(w_{00}^{2}w_{01}^{4}-w_{01}^{2}w_{11}^{4}\right)^{3}\left(a_{111}-w_{11}^{6}\right)}{a_{011}^{3}\left(w_{01}^{6}-w_{11}^{6}\right)}\\
\frac{\left(w_{00}^{4}w_{01}^{2}-w_{01}^{4}w_{11}^{2}\right)^{3}\left(w_{00}^{6}-a_{000}\right)}{a_{001}^{3}\left(w_{00}^{6}-w_{01}^{6}\right)} & = & \frac{\left(w_{00}^{2}w_{01}^{4}-w_{01}^{2}w_{11}^{4}\right)^{3}\left(w_{01}^{6}-a_{111}\right)}{a_{011}^{3}\left(w_{01}^{6}-w_{11}^{6}\right)}
\end{array}\end{cases}
\end{equation}
from which we deduce that the characteristic polynomial for $2\times2\times2$
hypermatrices is given by 
\begin{equation}
\left(w_{00}^{6}w_{11}^{6}-w_{01}^{12}\right)+w_{01}^{6}\left(a_{000}+a_{111}\right)-\left(a_{111}w_{00}^{6}+a_{000}w_{11}^{6}\right).
\end{equation}
Incidentally it immediately follows that characteristic polynomials
of direct sums of $2\times2\times2$ matrices is determined by the
derivation described above.

In order to derive the elimination ideal $I_{\mathbf{Q}}$ using the
hypermatrix formulation of the resolution identity we will consider
the sequence of hypermartrices defined as follows 
\begin{equation}
\mathbf{U}_{0}=\boldsymbol{\Delta},\qquad\mathbf{U}_{k+1}=\circ_{\mathbf{U}_{k}}\left(\mathbf{Q},\mathbf{Q}^{T^{2}},\mathbf{Q}^{T}\right)
\end{equation}
where $\boldsymbol{\Delta}$ denotes the Kronecker delta and the ternary
product determining $\mathbf{U}_{k+1}$ corresponds to the hypermatrix
product with background hypermatrix $\mathbf{U}_{k}$ as introduced
in \cite{GER}. The $k$-th term of the recurrence yields $n+2{n \choose 2}+2{n \choose 3}$
constraints, furthermore we know that for the purposes of elimination,
the number of variables being considered equals $n\left(n+2{n \choose 2}+2{n \choose 3}\right)$,
it therefore follows that it is enough to compute a sequence of length
$n$. Using the constraints 
\begin{equation}
\left\{ \mathbf{A}=\circ_{\mathbf{U}_{k}}\left(\circ\left(\mathbf{Q},\mathbf{D},\mathbf{D}^{T}\right),\,\circ\left(\mathbf{Q},\mathbf{D},\mathbf{D}^{T}\right)^{T^{2}},\,\circ\left(\mathbf{Q},\mathbf{D},\mathbf{D}^{T}\right)^{T}\right)\right\} _{0\le k<n}
\end{equation}
in conjunction with Cramer's rule we express the monomials in the
entries of $\mathbf{D}$ as rational function in the entries of $\mathbf{Q}$,
just as we did for matrices. Finally the elimination ideal $I_{\mathbf{Q}}$
is determined by the constraints of the form 
\begin{equation}
\begin{cases}
\begin{array}{ccc}
\left(\lambda_{ki}^{4}\lambda_{kp}^{2}\right)^{3} & = & \left[\left(\lambda_{ki}^{2}\right)^{3}\right]^{2}\left[\left(\lambda_{kp}^{2}\right)^{3}\right]\\
\left(\lambda_{ki}^{2}\lambda_{kj}^{2}\lambda_{kp}^{2}\right)^{3} & = & \left(\lambda_{ki}^{2}\right)^{3}\left(\lambda_{kj}^{2}\right)^{3}\left(\lambda_{kp}^{2}\right)^{3}
\end{array}\end{cases},\forall\,\left(i,j,p\right)\in\left\{ \underbrace{{l \choose 2}}_{\left(i,\, j,\, j\right)}\cup\underbrace{{l \choose 2}}_{\left(i,\, i,\, j\right)}\cup\underbrace{{l \choose 3}}_{\left(i,\, j,\, k\right)}\cup\underbrace{{l \choose 3}}_{\left(j,\, i,\, k\right)}\right\} 
\end{equation}
which determines the elimination ideal $I_{\mathbf{Q}}$. Incidentally
the necessary and sufficient condition for the existence of a spectral
decomposition for a given symmetric hypermatrix $\mathbf{A}$ is the
fact that the elimination ideal $I_{\mathbf{Q}}$ is non-trivial.

\subsection{Spectral bound.}

We recall that for a positive definite symmetric matrix $\mathbf{A}$
whose spectral decomposition is expressed by 
\begin{equation}
\begin{cases}
\begin{array}{ccc}
\mathbf{A} & = & \left(\mathbf{Q}\cdot\sqrt{\mathbf{D}}\right)\cdot\left(\mathbf{Q}\cdot\sqrt{\mathbf{D}}\right)^{T}\\
\left[\mathbf{Q}\cdot\mathbf{Q}^{T}\right]_{i,j} & = & \begin{cases}
\begin{array}{cc}
1 & \mbox{if }i=j\\
0 & \mbox{otherwise }
\end{array} & \forall\:0\le i,j<n\end{cases}\\
\mathbf{D}^{\star^{2}} & = & \mathbf{D}^{T}\cdot\mathbf{D}
\end{array}\end{cases},
\end{equation}
we have that 
\begin{equation}
\forall,\:\mathbf{x},\,\mathbf{y}\in\mathbb{C}^{n}\quad\left\langle \mathbf{x},\,\mathbf{y}\right\rangle _{\mathbf{A}}:=\left(\sum_{0\le i,j<n}a_{ij}x_{i}y_{j}\right)=\sum_{0\le k<n}\left\langle \sqrt{\lambda_{k}}\mathbf{x},\,\sqrt{\lambda_{k}}\mathbf{y}\right\rangle _{\mathbf{q}_{k}\otimes\mathbf{q}_{k}}
\end{equation}
it therefore follows from the resolution of identity that 
\begin{equation}
\forall,\:\mathbf{x},\,\mathbf{y}\in\mathbb{C}^{n}\quad\left\langle \mathbf{x},\,\mathbf{y}\right\rangle =\sum_{0\le k<n}\left\langle \mathbf{x},\,\mathbf{y}\right\rangle _{\mathbf{q}_{k}\otimes\mathbf{q}_{k}}
\end{equation}
if the vectors $\mathbf{x}$ and $\mathbf{y}$ are chosen such that
\[
\forall\:0\le k<n,\quad\left\langle \mathbf{x},\,\mathbf{y}\right\rangle _{\mathbf{q}_{k}\otimes\mathbf{q}_{k}}\ge0
\]
then the following spectral inequality holds 
\begin{equation}
\left\langle \sqrt{\lambda_{0}}\mathbf{x},\,\sqrt{\lambda_{0}}\mathbf{y}\right\rangle \le\left\langle \mathbf{x},\,\mathbf{y}\right\rangle _{\mathbf{A}}\le\left\langle \sqrt{\lambda_{n}}\mathbf{x},\,\sqrt{\lambda_{n}}\mathbf{y}\right\rangle 
\end{equation}
Similarly for some 3-hypermatrix $\mathbf{A}$ with entries symmetric
under cyclic permutation the corresponding spectral decomposition
is expressed by 
\begin{equation}
\begin{cases}
\begin{array}{ccc}
\mathbf{A} & = & \circ\left(\circ\left(\mathbf{Q},\mathbf{D},\mathbf{D}^{T}\right),\,\circ\left(\mathbf{Q},\mathbf{D},\mathbf{D}^{T}\right)^{T^{2}},\,\circ\left(\mathbf{Q},\mathbf{D},\mathbf{D}^{T}\right)^{T}\right)\\
\left[\circ\left(\mathbf{Q},\,\mathbf{Q}^{T^{2}},\,\mathbf{Q}^{T}\right)\right]_{i,j,k} & = & \begin{cases}
\begin{array}{cc}
1 & \mbox{if }i=j=k\\
0 & \mbox{otherwise }
\end{array} & \forall\:0\le i,j,k<n\end{cases}\\
\mathbf{D}^{\star^{3}} & = & \circ\left(\mathbf{D}^{T},\,\mathbf{D}^{T^{2}},\,\mathbf{D}\right)
\end{array}\end{cases}.
\end{equation}
We consider the very particular case where the scaling entries of
$\mathbf{D}$ are such that $\forall\,0\le j_{0}<j_{1}<n,\quad0\le d_{j_{0}}\left(i\right)\le d_{j_{1}}\left(i\right),$
we have that 

\[
\forall\:\mathbf{x},\,\mathbf{y},\,\mathbf{z}\in\mathbb{C}^{n},
\]
\begin{equation}
\left\langle \mathbf{x},\,\mathbf{y},\,\mathbf{z}\right\rangle _{\mathbf{A}}=\sum_{0\le k<n}\left\langle \left(\mathbf{d}_{k}\star\mathbf{x}\right),\,\left(\mathbf{d}_{k}\star\mathbf{y}\right),\,\left(\mathbf{d}_{k}\star\mathbf{z}\right)\right\rangle _{\otimes\left(Q_{k},Q_{k},Q_{k}\right)}
\end{equation}
$Q_{k}$ denote the $k$-th eigematrix of $\mathbf{A}$, $\otimes\left(Q_{k},Q_{k},Q_{k}\right)$
denote the matrix outer product as defined in \cite{GER} and $\mathbf{d}_{k}\star\mathbf{x}$
denotes the Hadamard product of the vectors $\mathbf{d}_{k}$ and
$\mathbf{x}$. Furthermore the hypermatrix resolution of identity
associated with the orthogonal hypermatrix $\mathbf{Q}$ is expressed
by 
\begin{equation}
\forall\:\mathbf{x},\,\mathbf{y},\,\mathbf{z}\in\mathbb{C}^{n},\quad\left\langle \mathbf{x},\,\mathbf{y},\,\mathbf{z}\right\rangle =\sum_{0\le j<n}\left\langle \mathbf{x},\,\mathbf{y},\,\mathbf{z}\right\rangle _{\otimes\left(Q_{j},Q_{j},Q_{j}\right)}
\end{equation}
in particular, if the vectors $\mathbf{x}$, $\mathbf{y}$ and $\mathbf{z}$
are chosen such that
\[
\forall\:0\le k<n,\quad\left\langle \left(\mathbf{d}_{k}\star\mathbf{x}\right),\,\left(\mathbf{d}_{k}\star\mathbf{y}\right),\,\left(\mathbf{d}_{k}\star\mathbf{z}\right)\right\rangle _{\otimes\left(Q_{k},Q_{k},Q_{k}\right)}\ge0
\]
then the following spectral inequality holds 
\begin{equation}
\left\langle \left(\mathbf{d}_{0}\star\mathbf{x}\right),\,\left(\mathbf{d}_{0}\star\mathbf{y}\right),\,\left(\mathbf{d}_{0}\star\mathbf{z}\right)\right\rangle \le\left\langle \mathbf{x},\,\mathbf{y},\,\mathbf{z}\right\rangle _{\mathbf{A}}\le\left\langle \left(\mathbf{d}_{n}\star\mathbf{x}\right),\,\left(\mathbf{d}_{n}\star\mathbf{y}\right),\,\left(\mathbf{d}_{n}\star\mathbf{z}\right)\right\rangle 
\end{equation}
which generalizes the matrix spectral inequality.

\subsection{3-hypermatrix SVD}

We now describe a natural generalization of symmetrization approach
to hypermatrix SVD. We start with some arbitrary $n\times n\times n$
hypermatrix $\mathbf{A}$, and deduce at most $3$ symmetric 3-hypermatrices
respectively given by $\circ\left(\mathbf{A},\,\mathbf{A}^{T^{2}},\,\mathbf{A}^{T}\right)$,
$\circ\left(\mathbf{A}^{T},\,\mathbf{A},\,\mathbf{A}^{T^{2}}\right)$
and $\circ\left(\mathbf{A}^{T^{2}},\,\mathbf{A}^{T},\,\mathbf{A}\right)$.
Furthermore as suggested by the spectral decomposition of 3-hypermatrices
which are symmetric under cyclic permutations of their indices we
are led to consider the following decomposition expressions associated
with each symmetric hypermatrices. 
\begin{equation}
\begin{cases}
\begin{array}{ccc}
\circ\left(\mathbf{A},\,\mathbf{A}^{T^{2}},\,\mathbf{A}^{T}\right) & = & \circ\left(\circ\left(\mathbf{Q},\mathbf{D},\mathbf{D}^{T}\right),\,\circ\left(\mathbf{Q},\mathbf{D},\mathbf{D}^{T}\right)^{T^{2}},\,\circ\left(\mathbf{Q},\mathbf{D},\mathbf{D}^{T}\right)^{T}\right)\\
\left[\circ\left(\mathbf{Q},\,\mathbf{Q}^{T^{2}},\,\mathbf{Q}^{T}\right)\right]_{i,j,k} & = & \begin{cases}
\begin{array}{cc}
1 & \mbox{if }i=j=k\\
0 & \mbox{otherwise }
\end{array} & \forall\:0\le i,j,k<n\end{cases}\\
\mathbf{D}^{\star^{3}} & = & \circ\left(\mathbf{D}^{T},\,\mathbf{D}^{T^{2}},\,\mathbf{D}\right)
\end{array}\end{cases}
\end{equation}
\begin{equation}
\begin{cases}
\begin{array}{ccc}
\circ\left(\mathbf{A}^{T},\,\mathbf{A},\,\mathbf{A}^{T^{2}}\right) & = & \circ\left(\circ\left(\mathbf{E},\mathbf{U},\mathbf{E}^{T^{2}}\right)^{T},\,\circ\left(\mathbf{E},\mathbf{U},\mathbf{E}^{T^{2}}\right),\,\circ\left(\mathbf{E},\mathbf{U},\mathbf{E}^{T^{2}}\right)^{T^{2}}\right)\\
\left[\circ\left(\mathbf{U}^{T},\,\mathbf{U},\,\mathbf{U}^{T^{2}}\right)\right]_{i,j,k} & = & \begin{cases}
\begin{array}{cc}
1 & \mbox{if }i=j=k\\
0 & \mbox{otherwise }
\end{array} & \forall\:0\le i,j,k<n\end{cases}\\
\mathbf{E}^{\star^{3}} & = & \circ\left(\mathbf{E}^{T},\,\mathbf{E}^{T^{2}},\,\mathbf{E}\right)
\end{array}\end{cases}
\end{equation}
\begin{equation}
\begin{cases}
\begin{array}{ccc}
\circ\left(\mathbf{A}^{T^{2}},\,\mathbf{A}^{T},\,\mathbf{A}\right) & = & \circ\left(\circ\left(\mathbf{F}^{T},\mathbf{F}^{T^{2}},\mathbf{V}\right)^{T^{2}},\,\circ\left(\mathbf{F}^{T},\mathbf{F}^{T^{2}},\mathbf{V}\right)^{T},\,\circ\left(\mathbf{F}^{T},\mathbf{F}^{T^{2}},\mathbf{V}\right)\right)\\
\left[\circ\left(\mathbf{V}^{T^{2}},\,\mathbf{V}^{T},\,\mathbf{V}\right)\right]_{i,j,k} & = & \begin{cases}
\begin{array}{cc}
1 & \mbox{if }i=j=k\\
0 & \mbox{otherwise }
\end{array} & \forall\:0\le i,j,k<n\end{cases}\\
\mathbf{F}^{\star^{3}} & = & \circ\left(\mathbf{F}^{T},\,\mathbf{F}^{T^{2}},\,\mathbf{F}\right)
\end{array}\end{cases}
\end{equation}
Incidentally the framework for the symmetrization approach results
from a desired product of the form 
\begin{equation}
\circ\left(\circ\left(\mathbf{Q},\mathbf{D},\mathbf{D}^{T}\right),\,\circ\left(\mathbf{E},\mathbf{U},\mathbf{E}^{T^{2}}\right),\,\circ\left(\mathbf{F}^{T},\mathbf{F}^{T^{2}},\mathbf{V}\right)\right).
\end{equation}
Let 
\[
\begin{array}{ccc}
\widetilde{\mathbf{Q}} & = & \circ\left(\mathbf{Q},\mathbf{D},\mathbf{D}^{T}\right)\\
\widetilde{\mathbf{E}} & = & \circ\left(\mathbf{E},\mathbf{U},\mathbf{E}^{T^{2}}\right)\\
\widetilde{\mathbf{F}} & = & \circ\left(\mathbf{F}^{T},\mathbf{F}^{T^{2}},\mathbf{V}\right)
\end{array}
\]
Ideally we would want to have 
\begin{equation}
\mathbf{A}=\circ\left(\widetilde{\mathbf{Q}},\,\widetilde{\mathbf{E}},\,\widetilde{\mathbf{F}}\right)
\end{equation}
but for most purposes we would be permit $n$ additional parameters
$\left\{ \alpha_{t}\right\} _{0\le t<n}$, which are to be solve in
the least square sense so as to yield an approximation of $\mathbf{A}$
\begin{equation}
\mathbf{A}\approx\sum_{0\le k<n}\alpha_{k}\quad\otimes\left(\widetilde{Q}_{k},\,\widetilde{E}_{k},\,\widetilde{T}_{k}\right)
\end{equation}

\section{General Hypermatrix Spectral Decomposition}

Just as matrices which are not symmetric admit a spectral decomposition,
3-hypermatrices which are not symmetric under cyclic permutation of
their indices also admit a spectral decomposition. We shall discuss
here the spectral decomposition of non-symmetric 3-hypermatrices.
We show here how to extend to arbitrary hypermatrices, the derivation
of the elimination ideals. For simplicity let us start with the matrix
case. We recall that matrix spectral constraint for symmetric real
matrices are expressed by 
\begin{equation}
\begin{cases}
\begin{array}{ccc}
\mathbf{A} & = & \left(\mathbf{U}\cdot\mathbf{D}\right)\cdot\left(\mathbf{V}\cdot\mathbf{D}\right)^{T}\\
\left[\mathbf{U}\cdot\mathbf{V}^{T}\right]_{i,j} & = & \begin{cases}
\begin{array}{cc}
1 & \mbox{if }i=j\\
0 & \mbox{otherwise }
\end{array} & \forall\:0\le i,j<n\end{cases}\\
\mathbf{D}^{\star^{2}} & = & \mathbf{D}^{T}\cdot\mathbf{D}
\end{array}\end{cases}.
\end{equation}
We describe here the computation of generators for the elimination
ideals
\begin{equation}
I_{\mathbf{U},\mathbf{V}}=\left\rangle \mathbf{A}-\left(\mathbf{U}\cdot\mathbf{D}\right)\cdot\left(\mathbf{V}\cdot\mathbf{D}\right)^{T},\,\mathbf{U}\cdot\mathbf{V}^{T}-\boldsymbol{\Delta},\,\mathbf{D}^{\star^{2}}-\mathbf{D}^{T}\cdot\mathbf{D}\right\langle \cap\mathbb{C}\left[\mathbf{U},\,\mathbf{V}\right]
\end{equation}
where $\left\rangle \mathbf{A}-\left(\mathbf{U}\cdot\mathbf{D}\right)\cdot\left(\mathbf{V}\cdot\mathbf{D}\right)^{T},\,\mathbf{U}\cdot\mathbf{V}^{T}-\boldsymbol{\Delta},\,\mathbf{D}^{\star^{2}}-\mathbf{D}^{T}\cdot\mathbf{D}\right\langle $
denotes the ideal generated by the matrix spectral constraints. Without
any loss of generality we may write the spectral constraints as 
\begin{equation}
\begin{cases}
\begin{array}{ccc}
a_{ij} & = & \left\langle \left(\mathbf{u}_{i}\star\boldsymbol{\lambda}\right),\left(\boldsymbol{\gamma}\star\mathbf{v}_{j}\right)\right\rangle \\
\delta_{ij} & = & \left\langle \mathbf{u}_{i},\mathbf{v}_{j}\right\rangle 
\end{array} & \forall\:0\le i,j<n\end{cases}
\end{equation}
Using the resolution of identity, we reduce the spectral constraints
to the following $n^{2}$ constraints of the form 
\begin{equation}
a_{ij}=\sum_{0\le k<n}\left\langle \left(\mathbf{u}_{i}\star\boldsymbol{\lambda}\right),\left(\boldsymbol{\gamma}\star\mathbf{v}_{j}\right)\right\rangle _{\mathbf{u}_{k}\cdot\mathbf{v}_{k}^{T}}
\end{equation}
which may more conveniently be rewritten as 
\begin{equation}
a_{ij}=\sum_{0\le k<n}\left\langle \left(\boldsymbol{\lambda}\cdot\boldsymbol{\gamma}^{T}\right),\,\left(\mathbf{u}_{k}\cdot\mathbf{v}_{k}^{T}\right)\star\left(\mathbf{u}_{i}\cdot\mathbf{v}_{j}^{T}\right)\right\rangle 
\end{equation}
which we think of as linear system of $n^{2}$ constraints in linear
in the $n^{2}$ variables $\left\{ \lambda_{i}\gamma_{j}\right\} _{0\le i,j<n}$.

Similarly for general hypermatrices the constraints is expressed by
\begin{equation}
\begin{cases}
\begin{array}{ccc}
\mathbf{A} & = & \circ\left(\circ\left(\mathbf{Q},\mathbf{D}_{0},\mathbf{D}_{0}^{T}\right),\,\circ\left(\mathbf{D}_{1},\mathbf{U},\mathbf{D}_{1}^{T^{2}}\right),\,\circ\left(\mathbf{D}_{2}^{T},\mathbf{D}_{2}^{T^{2}},\mathbf{V}\right)\right)\\
\left[\circ\left(\mathbf{Q},\,\mathbf{U},\,\mathbf{V}\right)\right]_{i,j,k} & = & \begin{cases}
\begin{array}{cc}
1 & \mbox{if }i=j=k\\
0 & \mbox{otherwise }
\end{array} & \forall\:0\le i,j,k<n\end{cases}\\
\mathbf{D}_{l}^{\star^{3}} & = & \circ\left(\mathbf{D}_{l}^{T},\,\mathbf{D}_{l}^{T^{2}},\,\mathbf{D}_{l}\right)\quad0\le l<3
\end{array}\end{cases}
\end{equation}
hence using the sequence 
\begin{equation}
\mathbf{G}_{0}=\boldsymbol{\Delta},\qquad\mathbf{G}_{k+1}=\circ_{\mathbf{G}_{k}}\left(\mathbf{Q},\mathbf{U},\mathbf{V}\right)
\end{equation}
Using the constraints 
\begin{equation}
\left\{ \mathbf{A}=\circ_{\mathbf{G}_{k}}\left(\circ\left(\mathbf{Q},\mathbf{D}_{0},\mathbf{D}_{0}^{T}\right),\,\circ\left(\mathbf{D}_{1},\mathbf{U},\mathbf{D}_{1}^{T^{2}}\right),\,\circ\left(\mathbf{D}_{2}^{T},\mathbf{D}_{2}^{T^{2}},\mathbf{V}\right)\right)\right\} _{0\le k<n}
\end{equation}
in conjunction with Cramer's rule just as we did for symmetric hypermatrices
we compute the elimination ideal for hypermatrices and deduce from
it a criteria for the existence of a spectral decomposition.

\section*{Acknowledgments}

This material is based upon work supported by the National Science
Foundation under agreements Princeton University Prime Award No. CCF-0832797
and Sub-contract No. 00001583. The author would like to thank the
IAS for providing excellent working conditions. The author is also
grateful to Vladimir Retakh, Avi Wigderson, Noga Alon for insightful
comments while preparing this manuscript.

\end{document}